\documentclass[journal, final, letterpaper, twocolumn, twoside]{IEEEtran}

\usepackage{cite} %\usepackage[noadjust]{cite}

\usepackage[cmex10]{amsmath}
\usepackage{amssymb, amsthm, eucal}
\usepackage{array}

\usepackage{graphicx, epsfig}
\usepackage[tight, footnotesize]{subfigure}

\newtheorem{proposition}{Proposition}
\newtheorem{example}{Example}
\newtheorem{assumption}{Assumption}

\def\g{\gamma}
\def\re{\Re}
\def\rn{\Re^n}

\def\pttj[#1]{\frac{\partial #1}{\partial \theta_j}}
\def\tpttj[#1]{\tfrac{\partial #1}{\partial \theta_j}}
\def\tpttk[#1]{\tfrac{\partial #1}{\partial \theta_k}}
\def\tpttjk[#1]{\tfrac{\partial^2 #1}{\partial \theta_j \partial \theta_k}}

\def \b {\beta}
\def\d{\delta}

\def \l {\lambda}
\def \m {\mu}

\def\smskip{\smallskip}
\def\pn{\par\noindent}
\def\lf{\left}
\def\ri{\right}
\def\tl{\tilde}
\def\ol{\overline}

\def\old#1{}

\def\eqref#1{(\ref{#1})}

\def\jstar{J^*}
\def\p{\pi}
\def\P{\Pi}
\def\qstar{Q^{\raise0.04pt\hbox{\sevenpoint *}} }

\def\en{{\cal E}_+^n}

\begin{document}

\title{Affine Monotonic and Risk-Sensitive Models in Dynamic Programming}

\author{Dimitri~P.~Bertsekas%
\thanks{}%
\thanks{D.\ P.\ Bertsekas is with the Laboratory for Information and Decision
Systems (LIDS), M.I.T.  Email: dimitrib@mit.edu.
Address: 77 Massachusetts  Ave.,\ Rm 32-D660, Cambridge, MA 02139, USA.}%  
}

\markboth{Report~LIDS-3204,~June 2016~ (Revised,~November~2017)}%
{Bertsekas: Affine Monotonic and Risk-Sensitive Models in Dynamic Programming}

 \maketitle

\begin{abstract}
In this paper we consider a broad class of infinite horizon discrete-time  optimal control models that involve a nonnegative cost function and an affine mapping in their dynamic programming equation. They include as special cases  several classical models such as  stochastic undiscounted nonnegative cost problems, stochastic multiplicative cost problems, and risk-sensitive problems with exponential cost. 
We focus on the case where the state space is finite and the control space has some compactness properties.  We assume that the affine mapping has a semicontractive character, whereby for some policies it is a contraction, while for others it is not. In one line of analysis, we impose assumptions guaranteeing that the noncontractive policies cannot be optimal. Under these assumptions, we prove strong results that resemble those for discounted Markovian decision problems, such as the uniqueness of solution of Bellman's equation, and the validity of forms of value and policy iteration. In the absence of these assumptions, the results are weaker and  unusual in character: the optimal cost function need not be a solution of Bellman's equation, and an optimal policy may not be found by value or policy iteration. Instead the optimal cost function over just the contractive policies  is the largest solution of Bellman's equation, and can be computed by a variety of algorithms.
\end{abstract}


\section{Introduction} \label{sec-intro}

We consider an  infinite horizon optimal control model, characterized by an affine and monotone abstract mapping that underlies the associated Bellman equation of dynamic programming (DP for short). This model was formulated with some analysis in the author's monograph [1] as a special case of abstract DP. In the present paper we will provide a deeper analysis  and more effective algorithms for the finite-state version of the model, under considerably weaker assumptions. 

To relate our analysis with the existing literature, we note that DP models specified by an abstract mapping defining the corresponding Bellman equation have a long history. Models where this mapping is a sup-norm contraction over the space of bounded cost functions were introduced by Denardo [2]; see also Denardo and Mitten  [3]. Their main area of application is discounted DP models of various types. Noncontractive models, where the abstract mapping is not a contraction of any kind, but is instead monotone, were considered by Bertsekas [4], [5] (see also Bertsekas and Shreve [6], Ch.\ 5). Among others, these models cover the important cases of positive and negative (reward) DP problems of Blackwell [7] and Strauch [8], respectively. Extensions of the analysis of [5] were given by Verdu and Poor [9], which considered additional structure that allows the development of backward and forward value iterations, and in the thesis by Szepesvari [10], [11], which introduced non-Markovian policies into the abstract DP framework.
The model of [5] was also used to develop asynchronous value iteration methods for abstract contractive and noncontractive DP models; see [12], 13]. Moreover, there have been extensions of the theory to asynchronous policy iteration algorithms and approximate DP by Bertsekas and Yu ([13], [14], [15]). 

A type of abstract DP model, called {\it semicontractive\/},  was introduced in the monograph [1]. In this model, the abstract DP mapping corresponding to some policies has a regularity/contraction-like property, but the mapping of others does not. A prominent example is the stochastic shortest path problem (SSP for short), a Markovian decision problem where we aim to drive the state of a finite-state Markov chain to a cost-free and absorbing termination state at minimum expected cost.  In SSP problems, the contractive policies are the so-called {\it proper} policies, which are the ones that lead to the termination state with probability 1. 
The SSP problem, originally introduced by Eaton and Zadeh [16], has been discussed under a variety of assumptions, in many sources, including the books [17], [18], [19], [20], [21], [22], [23], [24], and [12], where it is sometimes referred to by other names such as ``first passage problem" and ``transient programming problem." It has found a wide range of applications in regulation, path planning, robotics, and other contexts.

In this paper we focus on a different  subclass of semicontractive models, called {\it affine monotonic\/}, where the abstract mapping associated with a stationary policy is affine and maps nonnegative functions to nonnegative functions. These models include as special cases  stochastic undiscounted nonnegative cost problems (including SSP problems with nonnegative cost per stage), and multiplicative cost problems, such as problems with exponentiated cost. A key idea in our analysis is to use the notion of a {\it contractive policy} (one whose affine mapping involves a matrix with eigenvalues lying strictly within the unit circle). This notion is analogous to the one of a proper policy in SSP problems and is used in similar ways. 

Our analytical focus is on the validity and the uniqueness of solution of Bellman's equation, and the convergence of (possibly asynchronous) forms of value and policy iteration. Our results are analogous to those obtained for SSP problems by Bertsekas and Tsitsiklis [25], and Bertsekas and Yu [26].  As in the case of [25], under favorable assumptions where noncontractive policies cannot be optimal, we show that the optimal cost function is the unique solution of Bellman's equation, and we derive strong algorithmic results.  As in the case of [26], we consider more general situations where the optimal cost function need not be a solution of Bellman's equation, and  an optimal policy may not be found by value or policy iteration. To address such anomalies, we focus attention on the optimal cost function over just the contractive policies, and we show that it is the largest solution of Bellman's equation and that it is the natural limit of value iteration. 

However, there are some substantial differences from the analyses of [25] and [26]. The framework of the present paper is broader than SSP and  includes in particular multiplicative cost problems. Moreover, some of the assumptions are different and necessitate a different line of analysis; for example  there is no counterpart of the assumption that the optimal cost function is real-valued, which is fundamental in the analysis of [26]. As an indication, we note that deterministic shortest path problems with negative cost cycles can be readily treated within the framework of the present paper, but cannot be analyzed as SSP problems within the standard framework of [25] and the weaker framework of [26] because their optimal shortest path length is equal to $-\infty$ for some initial states. 

In this paper, we also pay special attention to  exponential cost problems, extending significantly some of the classical results of Denardo and Rothblum [27], and the more recent results of Patek [28]. Both of these papers impose assumptions that guarantee that the optimal cost function is the unique solution of Bellman's equation, whereas our assumptions are much weaker. The paper by  Denardo and Rothblum [27] also assumes a finite control space in order to bring to bear a line of analysis based on linear programming (see also the discussion of Section II).

The paper is organized as follows. In Section II we introduce the affine monotonic model, and we show that it contains as a special case multiplicative and exponential cost models. We also introduce contractive policies and related assumptions. In Section III we address the core analytical questions relating to Bellman's equation and its solution, and we obtain favorable results under the assumption that all noncontractive policies have infinite cost starting from some initial state. In Section IV we remove this latter assumption, and we show favorable results relating to a restricted  problem whereby we optimize over the contractive policies only. Algorithms such as value iteration, policy iteration, and linear programming are discussed somewhat briefly in this paper, since their analysis follows to a great extent established paths for semicontractive  abstract DP models [1].

Regarding notation, we denote by $\rn$ the standard Euclidean space of vectors $J=\big(J(1),\ldots,J(n)\big)$ with real-valued components, and we denote by $\re$ the real line.  We denote by $\rn_+$ the set of  vectors with nonnegative  real-valued components, $$\rn_+=\big\{J\mid 0\le J(i)<\infty,\ i=1,\ldots,n\big\},$$
and by ${\cal E}_+^n$ the set of
vectors with nonnegative extended real-valued components,
$${\cal E}_+^n=\big\{J\mid 0\le J(i)\le\infty,\ i=1,\ldots,n\big\}.$$
Inequalities with vectors are meant to be componentwise, i.e., $J\le J'$ means that $J(i)\le J'(i)$ for all $i$.
Similarly, in the absence of an explicit statement to the contrary, operations on vectors, such as $\lim$, $\limsup$, and $\inf$, are meant to be componentwise.

\section{Problem Formulation}

We consider a finite state space $X=\{1,\ldots,n\}$ and a (possibly infinite) control constraint set $U(i)$ for each state $i$. Let ${\cal M}$ denote the set of all functions $\m=\big(\m(1),\ldots,\m(n)\big)$ such that $\m(i)\in U(i)$ for each $i=1,\ldots,n$. By a {\it policy} we mean a sequence of the form $\pi=\{\mu_0,\mu_1,\ldots\}$, where $\mu_k\in{\cal M}$ for all $k=0,1,\ldots$. By a {\it stationary policy} we mean a policy of the form $\{\m,\m,\ldots\}$. For convenience we also refer to any $\m\in{\cal M}$ as a ``policy" and use it in place of the stationary policy $\{\m,\m,\ldots\}$, when confusion cannot arise. 

We introduce for each $\m\in{\cal M}$ a mapping $T_\m:{\cal E}_+^n\mapsto{\cal E}_+^n$ given by
\begin{equation} \label{eq-mumap}
T_\m J=b_\m+A_\m J,
\end{equation}
where $b_\m$ is a  vector of $\rn$ with components $b\big(i,\m(i)\big)$, $i=1,\ldots,n$, and $A_\m$ is an $n\times n$ matrix with components $A_{ij}\big(\m(i)\big)$, $i,j=1,\ldots,n$. We assume that
$b(i,u)$ and $A_{ij}(u)$ are nonnegative,
\begin{equation} \label{eq-affmonnoneg}
b(i,u)\ge0,\quad A_{ij}(u)\ge0,\quad \forall\ i,j=1,\ldots,n,\ u\in U(i).
\end{equation}
We define the mapping $T:{\cal E}_+^n\mapsto{\cal E}_+^n$, where for each $J\in {\cal E}_+^n$, $TJ$ is the vector of ${\cal E}_+^n$ with components
\begin{align} \label{eq-tmapaffmon}
(TJ)(i)=\inf_{u\in U(i)}&\lf[b(i,u)+\sum_{j=1}^nA_{ij}(u)J(j)\ri],\notag\\
&\qquad \qquad \qquad i=1,\ldots,n.
\end{align}
Note that since the value of the expression in braces on the right depends on $\m$ only through $\m(i)$, which is just restricted to be in $U(i)$,
we have
$$(TJ)(i)=\inf_{\m\in{\cal M}}(T_{\m} J)(i),\qquad  i=1,\ldots,n,$$
so that $(TJ)(i)\le (T_{\m} J)(i)$ for all $i$ and $\m\in{\cal M}$.

We now define a DP-like optimization problem that involves the mappings $T_\m$. We introduce a special vector $\bar J\in\rn_+$, and we define the cost function of a policy $\p=\{\m_0,\m_1,\ldots\}$ in terms of the composition of the mappings $T_{\m_k}$, $k=0,1,\ldots$, by
$$J_\p(i)=\limsup_{N\to\infty}\,(T_{\m_0}\cdots T_{\m_{N-1}} \bar J)(i),\qquad i=1,\ldots,n.$$
The cost function of a stationary policy $\m$, is written as
$$J_\m(i)=\limsup_{N\to\infty}\,(T_{\m}^N \bar J)(i),\qquad i=1,\ldots,n.$$
(We use $\limsup$ because we are not assured that the limit exists; our analysis and results remain essentially unchanged if $\limsup$ is replaced by $\liminf$.) In contractive abstract DP models, $T_\m$ is assumed to be a contraction for all $\m\in{\cal M}$, in which case $J_\m$ is the unique fixed point of $T_\m$ and does not depend on the choice of $\bar J$. Here we will not be making such an assumption, and the choice of $\bar J$ may affect profoundly the character of the problem. For example, in SSP and other additive cost Markovian decision problems $\bar J$ is the zero function, $\bar J(i)\equiv0$, while in multiplicative cost models $\bar J$ is the unit function, $\bar J(i)\equiv1$, as we will discuss shortly. Also in SSP problems $A_\m$ is a substochastic matrix for all $\m\in{\cal M}$,  while in other problems $A_\m$ can have components or row sums that are larger and as well as smaller than 1. 

We define the optimal cost function $\jstar$ by
$$\jstar(i)=\inf_{\p\in\P}J_\p(i),\qquad i=1,\ldots,n,$$
where $\P$ denotes the set of all policies.
We wish to find $\jstar$ and a  policy $\p^*\in\P$ that is optimal, i.e., $J_{\p^*}=\jstar$.  
The analysis of affine monotonic problems revolves around the equation $J=TJ$, or equivalently
\begin{equation} \label{eq-beaffmon}
J(i)=\inf_{u\in U(i)}\lf[b(i,u)+\sum_{j=1}^nA_{ij}(u)J(j)\ri],\quad  j=1,\ldots,n.
\end{equation}
This is the analog of the classical infinite horizon DP equation and it is referred to as {\it Bellman's equation\/}. We are interested in solutions of this equation (i.e., fixed points of $T$) within ${\cal E}_+^n$ and  within $\rn_+$. Usually in DP models one expects that $\jstar$ solves Bellman's equation, while optimal stationary policies can be obtained by minimization over $U(i)$ in its right-hand side. However, this is not true in general, as we will show in Section IV.

Affine monotonic models appear in several contexts. In particular, finite-state sequential stochastic control problems (including SSP problems) with nonnegative cost per stage (see, e.g., [12], Chapter 3, and Section IV) are special cases where $\bar J$ is the identically zero function [$\bar J(i)\equiv0$]. Also, discounted problems involving state and control-dependent discount factors (for example semi-Markov problems, cf.\ Section 1.4 of [12], or Chapter 11 of [29]) are special cases, with the discount factors being absorbed within the scalars $A_{ij}(u)$. In all of these cases, $A_\m$ is a substochastic matrix. There are also other special cases, where $A_\m$ is not substochastic. They correspond to interesting classes of practical problems, including SSP-type problems with a multiplicative or an exponential (rather than additive) cost function, which we proceed to discuss.

\subsection{Multiplicative and Exponential Cost SSP Problems}

We will describe a  type of  SSP problem, where the cost function of a policy accumulates over time multiplicatively, rather than additively, up to the termination state. The special case where the cost from a given state is the expected value of the exponential of the length of the path from the state up to termination was studied by Denardo and Rothblum [27], and Patek [28]. We are not aware of a study of the multiplicative cost version for problems where a cost-free and absorbing termination state plays a major role (the paper by Rothblum [30] deals with multiplicative cost problems but focuses on the average cost case).

Let us introduce in addition to the states $i=1,\ldots,n$, a cost-free and absorbing state $t$. There are probabilistic state transitions among the states $i=1,\ldots,n$, up to the first time a transition to state $t$ occurs, in which case the state transitions terminate. We denote by $p_{it}(u)$ and $p_{ij}(u)$ the probabilities of transition from $i$ to $t$ and to $j$ under $u$, respectively, so that
$$p_{it}(u)+\sum_{j=1}^np_{ij}(u)=1,\qquad i=1,\ldots,n,\ u\in U(i).$$
Next we introduce nonnegative scalars $h(i,u,t)$ and  $h(i,u,j)$,
$$h(i,u,t)\ge 0,\quad  h(i,u,j)\ge0,\quad \forall\ i,j=1,\ldots,n,\ u\in U(i),$$
and we consider the affine monotonic problem where the scalars $A_{ij}(u)$ and $b(i,u)$ are defined by
\begin{equation} \label{eq-hmultit}
A_{ij}(u)=p_{ij}(u)h(i,u,j),\qquad i,j=1,\ldots,n,\ u\in U(i),
\end{equation}
and 
\begin{equation} \label{eq-hmultio}
b(i,u)=p_{it}(u)h(i,u,t),\quad i=1,\ldots,n,\ u\in U(i),
\end{equation}
and the vector $\bar J$ is the unit vector,
$$\bar J(i)=1,\qquad i=1,\ldots,n.$$
The cost function of this problem has a multiplicative character as we show next.

Indeed, with the preceding definitions of $A_{ij}(u)$, $b(i,u)$, and $\bar J$, we will prove that the expression for the cost function of a policy $\p=\{\m_0,\m_1,\ldots\}$,
$$J_\p(x_0)=\limsup_{N\to\infty}\,(T_{\m_0}\cdots T_{\m_{N-1}} \bar J)(x_0),\qquad x_0=1,\ldots,n,$$
can be written in the multiplicative form
\begin{equation} \label{eq-infmulticost}
J_\p(x_0)=\limsup_{N\to\infty}\,E\lf\{\prod_{k=0}^{N-1}h\big(x_k,\m_k(x_k),x_{k+1}\big)\ri\},
\end{equation}
where:
\begin{itemize}
\item [(a)] $\{x_0,x_1,\ldots\}$ is the random state trajectory generated starting from $x_0$ and using $\p$.
\item [(b)] The expected value is with respect to the probability distribution of that trajectory.
\item [(c)] We use the notation 
$$h\big(x_k,\m_k(x_k),x_{k+1}\big)=1,\qquad \hbox{if }x_k=x_{k+1}=t,$$
(so that the multiplicative cost accumulation stops once the state reaches $t$). 
\end{itemize}
\pn Thus, we claim that {\it $J_\p(x_0)$ can be viewed as the expected value of cost accumulated multiplicatively, starting from $x_0$ up to reaching the termination state $t$ (or indefinitely accumulated multiplicatively, if $t$ is never reached)\/}.

To verify the  formula \eqref{eq-infmulticost} for $J_\p$, we use the definition 
$T_\m J=b_\m+A_\m J,$
to show by induction that for every $\p=\{\m_0,\m_1,\ldots\}$, we have
\begin{align} \label{eq-finhorcost}
T_{\m_0}\cdots T_{\m_{N-1}} \bar J&=A_{\m_0}\cdots A_{\m_{N-1}}\bar J\notag\\
&\qquad +b_{\m_0}+\sum_{k=1}^{N-1}A_{\m_0}\cdots A_{\m_{k-1}}b_{\m_k}.
\end{align}
We then interpret the $n$ components of each vector on the right as conditional expected values of the expression
\begin{equation} \label{eq-multicost}
\prod_{k=0}^{N-1}h\big(x_k,\m_k(x_k),x_{k+1}\big)
\end{equation}
multiplied with the appropriate conditional probability. In particular:
\begin{itemize}
\item [(a)] The $i$th component of the vector $A_{\m_0}\cdots A_{\m_{N-1}}\bar J$ in Eq.\ \eqref{eq-finhorcost} is the conditional expected value of the expression \eqref{eq-multicost}, given that $x_0=i$ and $x_N\ne t$, multiplied with the conditional probability that $x_N\ne t$, given that $x_0=i$.
\item [(b)] The $i$th component of the vector $b_{\m_0}$ in Eq.\ \eqref{eq-finhorcost}  is the conditional expected value of the expression \eqref{eq-multicost}, given that $x_0=i$ and $x_1=t$, multiplied with the conditional probability that $x_1= t$, given that $x_0=i$.
\item [(c)] The $i$th component of the vector $A_{\m_0}\cdots A_{\m_{k-1}}b_{\m_k}$ in Eq.\ \eqref{eq-finhorcost}  is the conditional expected value of the expression \eqref{eq-multicost}, given that $x_0=i$, $x_1,\ldots,x_{k-1}\ne t$, and $x_k=t$, multiplied with the conditional probability that $x_1,\ldots,$ $x_{k-1}\ne t$, and $x_k=t$, given that $x_0=i$.
\end{itemize}
\pn By adding these conditional probability expressions, we obtain the $i$th component of the unconditional expected value
$$E\lf\{\prod_{k=0}^{N-1}h\big(x_k,\m_k(x_k),x_{k+1}\big)\ri\},$$
thus verifying the formula \eqref{eq-infmulticost}.

A special case of multiplicative cost problem is the {\it risk-sensitive SSP problem with exponential cost function\/}, where for all $i=1,\ldots,n,$ and $u\in U(i)$,
\begin{equation} \label{eq-hexponspec}
h(i,u,j)=\hbox{exp}\big({g(i,u,j)}\big),\quad j=1,\ldots,n, t,
\end{equation}
and the function $g$ can take both positive and negative values. 
The Bellman equation for this problem is
\begin{align} \label{eq-hexpon}
J(i)=&\inf_{u\in U(i)}\Bigg[p_{it}(u)\hbox{exp}\big({g(i,u,t)}\big)\notag\\
&\ \ \ +\sum_{j=1}^np_{ij}(u)\hbox{exp}\big({g(i,u,j)}\big)J(j)\Bigg],\quad i=1,\ldots,n.
\end{align}
Based on Eq.\ \eqref{eq-infmulticost}, we have that $J_\p(x_0)$ is the limit superior of the expected value of the exponential of the $N$-step additive finite horizon cost up to termination, i.e., 
$$\sum_{k=0}^{\bar k}g\big(x_k,\m_k(x_k),x_{k+1}\big),$$
 where $\bar k$ is equal to the first index prior to $N-1$ such that $x_{\bar k+1}=t$, or is equal to $N-1$ if there is no such index. The use of the exponential  introduces risk aversion, by assigning a strictly convex increasing  penalty for large rather than small cost of a trajectory up to termination (and hence a preference for small variance of the additive cost up to termination). 

In the cases where $0\le g$ or $g\le 0$, we also have $\bar J\le T\bar J$ and $T\bar J\le \bar J$, respectively, corresponding to a monotone increasing and a monotone decreasing problem, in the terminology of [1]. Both of these problems admit a favorable analysis, highlighted by the fact that $\jstar$ is a fixed point of $T$ (see [1], Chapter 4). 

The case where $g$ can take both positive and negative values is more challenging, and is the focus of this paper. We will consider two cases, discussed in Sections III and IV of this paper, respectively. Under the assumptions of Section III, $\jstar$ is shown to be the unique fixed point of $T$ within $\rn_+$. Under the assumptions of Section IV,  it may happen that $\jstar$ is not a fixed point of $T$ (see Example \ref{examplebelcounterexp} that follows). Denardo and Rothblum [27] and Patek [28] consider only the more benign Section III case, for which $\jstar$ is a fixed point of $T$. Also, the approach of [27]  is very different from ours: it relies on linear programming ideas, and for this reason it requires a finite control constraint set and cannot be readily adapted to an infinite control space. The approach of [28] is closer to ours in that it also descends from the paper [25]. It allows for an infinite control space under a compactness assumption that is similar to our Assumption 2.2 of the next section, but it also requires that $g(i,u,j)>0$ for all $i,u,j$, so it deals only with a monotone increasing case where $T_\m \bar J\ge \bar J$ for all $\m\in{\cal M}$. 

The deterministic version of the exponential cost problem where for each $u\in U(i)$, only  one of the transition probabilities $p_{it}(u),p_{i1}(u),\ldots,p_{in}(u)$ is equal to 1 and all others are equal to 0, is mathematically equivalent to the classical deterministic shortest path problem (since minimizing the exponential of a deterministic expression is equivalent to minimizing that expression). For this problem a standard assumption is that there are no cycles that have negative total length to ensure that the shortest path length is finite. However, it is interesting that this assumption is not required in the present paper: when there are paths that travel perpetually around a negative length cycle we simply have $\jstar(i)=0$ for all states $i$ on the cycle, which is permissible within our context.

\subsection{Assumptions on Policies - Contractive Policies}

\pn  We now introduce a special property of policies which is central for the purposes of this paper. We say that a given stationary policy $\m$ is {\it contractive if $A_\m^N\to0$ as $N\to\infty$\/}. Equivalently, $\m$ is contractive if all the eigenvalues of $A_\m$ lie strictly within the unit circle. Otherwise, $\m$ is called {\it noncontractive\/}. It follows from a classical result that a policy $\m$ is contractive if and only if $T_\m$ is a contraction with respect to some norm. Because $A_\m\ge0$, a stronger assertion can be made: $\m$ is contractive if and only if $A_\m$ is a contraction with respect to some weighted sup-norm (see e.g., the discussion in [22], Ch.\ 2, Cor.\ 6.2, or [12], Section 1.5.1). In the special case of SSP problems with additive cost function (all matrices $A_\m$ are substochastic), the contractive policies coincide with the {\it proper} policies, i.e., the ones that lead to the termination state with probability 1, starting from every state. 

A particularly favorable situation for an SSP problem arises when all policies are proper, in which case all the mappings $T$ and $T_\m$ are contractions with respect to some common weighted sup norm. This result was shown in the paper by Veinott [31], where it was attributed to A.\ J.\ Hoffman. Alternative proofs of this contraction property are given in Bertsekas and Tsitsiklis, [22], p.\ 325 and  [32], Prop.\ 2.2, Tseng [33],  and Littman [34]. The proofs of [32] and [34] are essentially identical, and  easily generalize to the context of the present paper. In particular, it can be shown that 
if all policies are contractive,  all the mappings $T$ and $T_\m$ are contractions with respect to some common weighted sup norm. However, we will not prove or use this fact in this paper.

Let us derive an expression for the cost function of contractive and noncontractive policies.
By repeatedly applying the mapping $T$ to both sides of the equation $T_\m J=b_\m+A_\m J$, we have
$$T^N_\m J=A_\m^N  J+\sum_{k=0}^{N-1} A_\m^k b_\m,\qquad \forall\ J\in \rn,\ N=1,2,\ldots,$$
and hence
\begin{equation} \label{eq-infhorcost}
J_\m=\limsup_{N\to\infty}T_\m^N\bar J=\limsup_{N\to\infty}A_\m^N  \bar J+\sum_{k=0}^{\infty} A_\m^k b_\m.
\end{equation}
From these expressions,
it follows that if $\m$ is contractive, the initial function $\bar J$ in the definition of $J_\m$ 
does not matter, and we have
\begin{align*}
J_\m&=\limsup_{N\to\infty}T_\m^NJ\notag\\
&=\limsup_{N\to\infty}\sum_{k=0}^{N-1} A_\m^k b_\m,\quad   \forall\ \m\hbox{: contractive},\ J\in\rn.
\end{align*}
Moreover, since for a contractive $\m$, $T_\m$ is a contraction with respect to a weighted sup-norm, the $\lim\sup$ above can be replaced by $\lim$, so that
\begin{equation} \label{eq-contractiveexpr}
J_\m=\sum_{k=0}^{\infty} A_\m^k b_\m=(I-A_\m)^{-1}b_\m,\qquad \forall\ \m\hbox{: contractive}.
\end{equation}
Thus if  $\m$ is contractive, $J_\m$ is real-valued as well as nonnegative, i.e., $J_\m\in\rn_+$. If $\m$ is noncontractive, we have $J_\m\in\en$ and it is possible that 
$J_\m(i)=\infty$ for some states $i$. We will assume throughout the paper the following.

\begin{assumption}\label{assumptiontwoone}
There exists at least one contractive policy.
\end{assumption}

The analysis of finite-state SSP problems typically assumes that the control space is either finite, or satisfies a compactness and continuity condition introduced in [25]. The following is a similar condition, and will be in effect throughout the paper. 

\begin{assumption}[Compactness and Continuity]\label{assumptiontwotwo}
The control space $U$ is a metric space, and $p_{ij}(\cdot)$ and $b(i,\cdot)$ are continuous functions of $u$ over $U(i)$, for all $i$ and $j$. Moreover, for each state $i$, the sets
$$\lf\{u\in U(i)\ \Big|\ b(i,u)+\sum_{j=1}^nA_{ij}(u)J(j)\le \l\ri\}$$
are compact subsets of $U$ for all scalars $\l\in\re$ and $J\in \rn_+$.
\end{assumption}
 
 The preceding assumption is satisfied if the control space $U$ is finite. One way to see this is to simply identify each $u\in U$ with a distinct integer from the real line. Another interesting case where the assumption is satisfied is when for all $i$, $U(i)$ is a compact subset of the metric space $U$, and the functions $b(i,\cdot)$ and $A_{ij}(\cdot)$ are continuous functions of $u$ over $U(i)$.

An advantage of allowing $U(i)$ to be infinite and compact is that it makes possible the use of randomized policies for problems where there is a {\it finite} set of feasible actions at each state $i$, call it $C(i)$. We may then specify $U(i)$ to be the set of all probability distributions over $C(i)$, which is a compact subset of a Euclidean space. In this way, our results apply to finite-state and finite-action problems where randomization is allowed, and $\jstar$ is the optimal cost function over all randomized nonstationary policies. Note, however, that the optimal cost function may change when randomized policies are introduced in this way. Basically, for our purposes, optimization over nonrandomized and over randomized policies over finite action sets $C(i)$ are two different problems, both of which are interesting and can be addressed with the methodology of this paper. However, when the sets $C(i)$ are infinite, a different and mathematically more sophisticated framework is required in order to allow randomized policies. The reason is that randomized policies over the infinite action sets $C(i)$ must obey measurability restrictions, such as universal measurability; see Bertsekas and Shreve [6], and Yu and Bertsekas [35].

The  compactness and continuity part of the preceding assumption guarantees some important properties of the mapping $T$. These are summarized in the following proposition.

\begin{proposition}\label{prp-propcompactimpl}
Let Assumptions \ref{assumptiontwoone} and \ref{assumptiontwotwo} hold.
\begin{itemize}
\item [(a)] The set of $u\in U(i)$ that minimize the expression 
\begin{equation} \label{eq-minimexpr}
b(i,u)+\sum_{j=1}^nA_{ij}(u)J(j)
\end{equation}
is nonempty and compact.
\item [(b)] Let $J_0$ be the zero vector in $\rn$ [$J_0(i)\equiv0$]. The sequence $\{T^kJ_0\}$ is monotonically nondecreasing and converges to a limit $\tl J\in \rn_+$ that satisfies $\tl J\le \jstar$ and $\tl J=T \tl J$.
\end{itemize}
\end{proposition}

\proof (a) The set of $u\in U(i)$ that minimize the expression in Eq.\ \eqref{eq-minimexpr} is the intersection $\cap_{m=1}^\infty U_m$ of the nested sequence of sets 
\begin{align*}
U_m=\lf\{u\in U(i)\ \Big|\ b(i,u)+\sum_{j=1}^nA_{ij}(u)J(j)\le\l_m\ri\},\\
m=1,2,\ldots,
\end{align*}
where $\{\l_m\}$ is a monotonically decreasing  sequence such that
$$\l_m\downarrow \inf_{u\in U(i)}\lf[b(i,u)+\sum_{j=1}^nA_{ij}(u)J(j)\ri].$$
Each set $U_m$ is nonempty, and by Assumption \ref{assumptiontwotwo}, it is compact, so the intersection is nonempty and compact.

\smskip
\pn (b) By the nonnegativity of $b(i,u)$ and $A_{ij}(u)$, we have $J_0\le TJ_0$, which by the monotonicity of $T$ implies that  $\{T^kJ_0\}$ is monotonically nondecreasing to a limit $\tl J\in\en$, and we have 
\begin{equation} \label{eq-incrsect}
J_0\le TJ_0 \le \cdots\le T^kJ_0\le \cdots\le
\tl J.
\end{equation}
For all policies $\p=\{\m_0,\m_1,\ldots\}$, we have $T^kJ_0\le T^k\bar J\le T_{\m_0}\cdots T_{\m_{k-1}}\bar J$, so by taking the limit as $k\to\infty$, we obtain $\tl J\le J_\p$, and  by taking the infimum over $\p$, it follows that $\tl J\le \jstar$. By Assumption \ref{assumptiontwoone}, there exists at least one contractive policy $\m$, for which $J_\m$ is real-valued [cf.\ Eq.\ \eqref{eq-contractiveexpr}], so $\jstar\in\rn_+$. It follows that the sequence $\{T^kJ_0\}$ consists of vectors in $\rn_+$, while $\tl J\in\rn_+$.

By applying $T$ to both sides of Eq.\ \eqref{eq-incrsect}, we obtain 
\begin{align*}
(T^{k+1}J_0)(i) &= \inf_{u\in U(i)}
\lf[b(i,u)+\sum_{j=1}^nA_{ij}(u)(T^kJ_0)(j)\ri]\\
&\le
(T\tl J)(i),
\end{align*} 
and by taking the limit as $k\to\infty$, it
follows that $\tl J\le T\tl J.$ 
Assume to arrive at a contradiction that there exists a state $\tl i$ such that
\begin{equation} \label{eq-onetwethsf}
\tl J(\tl i)< (T\tl J)(\tilde i).
\end{equation}
Consider the sets 
$$U_k(\tilde i)= \lf\{ u\in U(\tilde i)\ 
\Big|\ b(\tl i,u) + \sum_{j=1}^nA_{\tl ij}(u)(T^kJ_0)(j)\le \tl J(\tl i)\ri\}$$  
for $k\ge0$. It follows from  Assumption \ref{assumptiontwotwo} and Eq.\ \eqref{eq-incrsect} that $\big\{U_k(\tilde i)\big\}$ is a nested sequence of compact sets. Let
also $u_k$ be a control attaining the minimum in 
$$\min_{u\in U(\tilde
i)}\lf[b(\tl i,u)+\sum_{j=1}^nA_{\tl ij}(u)(T^kJ_0)(j)\ri];$$ 
[such a control exists by part (a)].  From Eq.\ \eqref{eq-incrsect}, it follows that for all $m\ge k$,
\begin{align*}
b(\tilde i,u_m)&+
\sum_{j=1}^nA_{\tl ij}(u_m)(T^kJ_0)(j) \\
&\le b(\tilde i,u_m) +
\sum_{j=1}^nA_{\tl ij}(u_m)(T^mJ_0)(j)\\
&\le \tl J(\tilde i).
\end{align*}  
Therefore $\{
u_m\}_{m=k}^\infty \subset U_k(\tilde i)$, and since $U_k(\tilde i)$ is compact, all the limit points of $\{ u_m\}_{m=k}^\infty$ belong to
$U_k(\tilde i)$ and at least one such limit point exists.  Hence the
same is true of the limit points of the entire sequence $\{ u_m\}_{m=0}^\infty$.  It
follows that if $\tilde u$ is a limit point of $\{ u_m\}_{m=0}^\infty$ then $$\tilde
u\in \cap_{k=0}^\infty U_k(\tilde i).$$ 
This implies that for all $k\ge 0$ 
$$(T^{k+1}J_0)(\tilde i)\le  b(\tilde i,\tilde u) + \sum_{j=1}^nA_{\tl ij}(\tilde u)(T^kJ_0)(j)\le \tl J(\tilde i).$$ 
By taking the limit in this relation as $k\to\infty$, we obtain 
$$\tl J(\tilde i) = b(\tilde i,\tilde u) + \sum_{j=1}^nA_{\tl ij}(\tilde u)\tl J(j).$$ 
Since the right-hand side is greater than or equal
to
$(T\tl J)(\tilde i)$, Eq.\ \eqref{eq-onetwethsf} is contradicted, implying that $\tl J =T\tl J$. \qed

\section{Case of Infinite Cost Noncontractive Policies}

\pn We now turn to  questions relating to Bellman's equation, the convergence of value iteration (VI for short) and policy iteration (PI for short), as well  as conditions for optimality of a stationary policy. In this section we will use the following assumption, which parallels the central assumption of [25] for SSP problems. We will not need this assumption in Section IV.

\begin{assumption}[Infinite Cost Condition] \label{assumptionthreeone}
For every noncontractive policy $\m$, there is at least one state such that the corresponding component of the vector
$\sum_{k=0}^\infty A_\m^k b_\m$ is equal to $\infty$.
\end{assumption}
 
 Note that the preceding assumption guarantees that for every noncontractive policy $\m$, we have $J_\m(i)=\infty$ for at least one state $i$ [cf.\ Eq.\ \eqref{eq-infhorcost}]. The reverse is not true, however: $J_\m(i)=\infty$ does not imply that the $i$th component of $\sum_{k=0}^\infty A_\m^k b_\m$ is equal to $\infty$, since there is the possibility that $A_\m^N\bar J$ may become unbounded as $N\to\infty$ [cf.\ Eq.\ \eqref{eq-infhorcost}]. 

Under Assumptions \ref{assumptiontwoone}, \ref{assumptiontwotwo}, and \ref{assumptionthreeone}, we will now derive results that closely parallel the standard results of [25] for additive cost SSP problems.
We have the following characterization of contractive policies.

\begin{proposition}[Properties of Contractive Policies]\label{prp-proptto}
Let Assumption \ref{assumptionthreeone} hold.
\begin{itemize}
\item [(a)] For a contractive policy  $\m$, the associated
cost vector $J_\m$ satisfies
$$\lim_{k\to\infty}(T_\m^k J)(i)=J_\m(i),\qquad
i=1,\ldots,n,$$ 
for every vector $J\in\rn$. Furthermore, we have 
$J_\m=T_\m J_\m,$
 and $J_\m$
is the unique solution of this equation within $\rn$.
\item [(b)] A stationary policy  $\m$ is contractive if and only if it satisfies $J\ge
T_\m J$
for some vector $J\in\rn_+$.
\end{itemize}
\end{proposition}

\proof (a) Follows from Eqs.\ \eqref{eq-finhorcost} and \eqref{eq-contractiveexpr}, and by writing the equation $J_\m=T_\m J_\m$ as $(I-A_\m) J_\m=b_\m$.
\smskip

\pn(b) If $\m$ is contractive, by part (a) we have $J\ge T_\m J$ for $J=J_\m\in \rn_+$. Conversely, let $J$ be a vector in $\rn_+$ with $J\ge T_\m J$. Then the monotonicity of $T_\m$ and Eq.\ \eqref{eq-finhorcost} imply that for all $N$ we have
 $$J\ge T_\m^N J=A_\m^NJ+\sum_{k=0}^{N-1}A_\m^k b_\m\ge \sum_{k=0}^{N-1}A_\m^k b_\m\ge0.$$ 
It follows that the vector $\sum_{k=0}^{\infty}A_\m^k b_\m$ is real-valued so that, by Assumption \ref{assumptionthreeone}, $\m$ cannot be noncontractive.  
\qed
\smskip

The following proposition is our main result under Assumption \ref{assumptionthreeone}. It parallels Prop.\ 3 of [25] (see also Section 3.2 of [12]). In addition to the fixed point property of $\jstar$ and the convergence of the VI sequence $\{T^kJ\}$ to $\jstar$ starting from any $J\in\rn_+$, it shows the validity of the PI algorithm. The latter algorithm generates a sequence $\{\m^k\}$ starting from any contractive policy $\m^0$. Its typical iteration consists of a computation of $J_{\m^k}$ using the policy evaluation equation $J_{\m^k}=T_{\m^k}J_{\m^k}$, followed by the policy improvement operation $T_{\m^{k+1}}J_{\m^k}=TJ_{\m^k}$.

\begin{proposition}[Bellman's Equation,  Policy Iteration, Value Iteration, and Optimality Conditions]\label{prp-propttt}
Let Assumptions \ref{assumptiontwoone}, \ref{assumptiontwotwo}, and \ref{assumptionthreeone}\ hold.
\begin{itemize}
\item [(a)] The optimal cost vector $\jstar $ satisfies the Bellman equation $J =T J$. 
Moreover, $\jstar $ is the unique solution of this equation within $\rn_+$. 
\item [(b)] Starting with any contractive policy $\m^0$, the sequence $\{\m^k\}$ generated by the PI algorithm consists of contractive policies, and any limit point of this sequence is a contractive optimal  policy.
\item [(c)] We have $$\lim_{k\to\infty}(T^kJ)(i)=\jstar (i),\qquad i=1,\ldots,n,$$ for every
vector $J\in\rn_+$.
\item [(d)] A stationary policy  $\m$ is optimal if and only if 
 $T_\m \jstar =T\jstar .$
\item [(e)] For a vector $J\in \rn_+$, if $J\le TJ$ then $J\le \jstar $, and if $J\ge TJ$ then $J\ge \jstar $.\end{itemize}
\end{proposition}

\proof (a), (b) From Prop.\ \ref{prp-propcompactimpl}(b), $T$ has as fixed point the vector $\tl J\in\rn_+$,  the limit of the sequence $\{T^kJ_0\}$, where $J_0$ is the identically zero vector [$J_0(i)\equiv0$]. We will show parts (a) and (b) simultaneously and in stages. First we will show that $\tl J$ is the unique fixed point of $T$ within $\rn_+$. Then we will show that the PI algorithm, starting from any contractive policy, generates in the limit a contractive policy $\ol \m$ such that $J_{\ol\m}=\tl J$. Finally we will show that $J_{\ol\m}=\jstar$.

Indeed, if $J$ and $J'$ are
two fixed points, then we select $\m$ and $\m'$ such that $J=TJ =T_\m J$ and
$J'=TJ' =T_{\m'}J'$; this is possible because of Prop.\ \ref{prp-propcompactimpl}(a). By Prop.\
\ref{prp-proptto}(b), we have that $\m$ and $\m'$ are contractive, and by Prop.\ \ref{prp-proptto}(a) we obtain $J=J_\m$ and
$J'=J_{\m'}$. We also have $J=T^kJ\le T_{\m'}^kJ$ for all $k\ge1$, and by Prop.\ \ref{prp-proptto}(a), it follows that 
$J\le\lim_{k\to\infty}T_{\m'}^k J=J_{\m'}=J'$. Similarly, $J'\le J$, showing  that   $J=J'$.
Thus $T$ has $\tl J$ as its unique fixed point within $\rn_+$.

We next turn to the PI algorithm. Let $\m$ be  a contractive policy (there
exists one by Assumption \ref{assumptiontwoone}). Choose $\m'$ such that $$T_{\m'} J_\m=TJ_\m.$$
 Then we have $J_\m=T_\m J_\m\ge T_{\m'}J_\m$. By Prop.\ \ref{prp-proptto}(b), $\m'$ is contractive, and using
the monotonicity of $T_{\m'}$ and Prop.\ \ref{prp-proptto}(a), we obtain
\begin{equation} \label{eq-ones}
J_\m\ge\lim_{k\to\infty}T_{\m'}^k J_\m=J_{\m'}.
\end{equation}
Continuing in the same manner, we construct a sequence $\{\m^k\}$ such that each $\m^k$ is contractive
and  
\begin{equation} \label{eq-onese}
J_{\m^k}\ge T_{\m^{k+1}}J_{\m^k} =TJ_{\m^k} \ge J_{\m^{k+1}},\qquad  k=0,1,\ldots
\end{equation}
The sequence $\{J_{\m^k}\}$ is real-valued, nonincreasing, and nonnegative so it converges to some $J_\infty\in\rn_+$.

We  claim that the sequence of vectors $\m^k=\big(\m^k(1),\ldots,\m^k(n)\big)$ has a limit point $\big(\ol \m(1),\ldots,\ol\m(n)\big)$, with $\ol\m$ being a feasible policy. Indeed, using Eq.\ \eqref{eq-onese} and the fact $J_\infty\le J_{\m^{k-1}}$, we have for all $k=1,2,\ldots,$
$$T_{\m^{k}}J_\infty\le T_{\m^{k}}J_{\m^{k-1}}=TJ_{\m^{k-1}}\le T_{\m^{k-1}}J_{\m^{k-1}}=J_{\m^{k-1}}\le J_{\m^0},$$
so $\m^{k}(i)$ belongs to the set
$$\hat U(i)=\lf\{u\in U(i)\ \Big|\  b(i,u)+\sum_{j=1}^nA_{ij}(u)J_\infty(j)\le J_{\m^0}(i)\ri\},$$
which is compact by Assumption \ref{assumptiontwotwo}. Hence the sequence $\{\m^k\}$ belongs to the compact set $\hat U(1)\times\cdots\times \hat U(n)$, and has a limit point $\ol\m$, which is a feasible policy.
In what follows, without loss of generality, we assume that the entire sequence $\{\m^k\}$ converges to $\ol\m.$ 

Since $J_{\m^k}\downarrow J_\infty\in \rn_+$ and $\m^k\to\ol \m$, by taking limit as $k\to\infty$ in Eq.\ \eqref{eq-onese}, and using the continuity part of Assumption \ref{assumptiontwotwo}, 
we obtain $J_\infty=T_{\ol\m} J_\infty$. It follows from Prop.\ \ref{prp-proptto}(b) that $\ol\m$ is contractive, and that $J_{\ol\m}$ is equal to $J_\infty$. To show that $J_{\ol\m}$ is a fixed point of $T$, we note that from the right side of Eq.\ \eqref{eq-onese}, we have for all policies $\m$, $T_\m J_{\m^k} \ge J_{\m^{k+1}}$, which by taking limit as $k\to\infty$ yields $T_\m J_{\ol\m}\ge J_{\ol\m}$. By taking minimum over $\m$, we obtain $T J_{\ol\m}\ge J_{\ol\m}$. Combining this with the relation $J_{\ol\m}=T_{\ol\m} J_{\ol\m}\ge TJ_{\ol\m}$, it follows that $J_{\ol\m}= TJ_{\ol\m}$. Thus $J_{\ol\m}$ is equal to the unique
fixed point $\tl J$ of $T$ within $\rn_+$.

We will now conclude the proof by showing that $J_{\ol\m}$ is equal to the optimal cost vector $\jstar $ (which also implies the optimality of the policy $\ol\m$, obtained from the  PI algorithm starting from a contractive policy). By Prop.\ \ref{prp-propcompactimpl}(b),  the sequence $T^kJ_0$ converges monotonically to  $\tl J$, which is equal to $J_{\ol\m}$. Also, for every policy $\p=\{\m_0,\m_1,\ldots\}$, we have
$$T^kJ_0\le T^k\bar J\le T_{\m_0}\cdots T_{\m_{k-1}}\bar J,\qquad k=0,1,\ldots,$$
and by taking the limit as $k\to\infty$, we obtain $J_{\ol \m}=\tl J=\lim_{k\to\infty}T^k J_0\le J_\p$ for all $\p$, showing that $J_{\ol\m}=\jstar$. Thus $\jstar$ is the unique fixed point of $T$ within $\rn_+$, and $\ol \m$ is an optimal policy.

\smskip
\pn (c) From the preceding proof, we have that 
$T^kJ_0\to \jstar$, which implies that 
\begin{equation} \label{eq-convbelow}
\lim_{k\to\infty}T^kJ=\jstar,\qquad \forall\ J\in\rn_+\hbox{ with }J\le \jstar.
\end{equation}
Also, for any $J\in\rn_+$ with $J\ge \jstar$, we have
$$T_{\ol\m}^kJ\ge T^kJ\ge T^k\jstar=\jstar=J_{\ol\m},$$
where $\ol\m$ is the contractive optimal policy obtained by PI in the proof of part (b). By taking the limit as $k\to\infty$ and using the fact $T_{\ol\m}^kJ\to J_{\ol\m}$ (which follows from the contractiveness of $\ol\m$), we obtain 
\begin{equation} \label{eq-convabove}
\lim_{k\to\infty}T^kJ=\jstar,\qquad \forall\ J\in\rn_+\hbox{ with }J\ge \jstar.
\end{equation}
Finally, given any $J\in\rn_+$, we have from Eqs.\ \eqref{eq-convbelow} and \eqref{eq-convabove},
$$\lim_{k\to\infty}T^k\big(\min\{J,\jstar\}\big)=\jstar,\quad \lim_{k\to\infty}T^k\big(\max\{J,\jstar\}\big)=\jstar,$$
and since $J$ lies between $\min\{J,\jstar\}$ and $\max\{J,\jstar\}$, it follows that $T^kJ\to\jstar$. 
\smskip

\pn(d) If  $\m$  is optimal, then $J_\m=\jstar $ and since by part (a) $\jstar$ is real-valued, $\m$ is
contractive. Therefore, by Prop.\ \ref{prp-proptto}(a), 
$$T_{\m}\jstar =T_{\m}J_{\m}=J_{\m}=\jstar =T\jstar  .$$
Conversely, if
$\jstar =T\jstar  =T_{\m}\jstar $, it follows from Prop.\ \ref{prp-proptto}(b) that   $\m$ is contractive, and by using
Prop.\ \ref{prp-proptto}(a), we obtain $\jstar =J_\m$. Therefore $\m$ is optimal.
The existence of an optimal policy follows from part (b).
\smskip

\pn(e) If $J\in \rn_+$ and $J\le TJ$, by repeatedly applying $T$ to both sides and using the monotonicity of $T$, we obtain
$J\le T^kJ$ for all  $k.$
Taking the limit as $k\to\infty$ and using the fact $T^kJ\to \jstar $ [cf.\ part (c)], we obtain $J\le \jstar $. The proof that $J\ge \jstar $ if $J\ge TJ$ is similar.  
 \qed

\subsection{Computational Methods}

 Proposition \ref{prp-propttt}(b) shows the validity of PI when starting from a contractive policy.  
This is similar to the case of additive cost SSP, where PI is known to converge starting from a proper policy (cf.\ the proof of Prop.\ 3 of [25]). 

There is also an asynchronous version of the PI algorithm proposed for discounted and SSP models by Bertsekas and Yu [14], [15], which does not require an initial contractive policy and admits an asynchronous implementation. This algorithm extends straightforwardly to the affine monotonic model of this paper  under Assumptions 2.1, 2.2, and 3.1 (see [1], Section 3.3.2, for a description of this extension to abstract DP models). 

Proposition \ref{prp-propttt}(c) establishes the validity of the VI algorithm that generates the sequence $\{T^kJ\}$, starting from any initial $J\in \rn_+$. An asynchronous version of this algorithm is also valid; see the discussion of Section 3.3.1 of [1].

Finally, Prop.\ \ref{prp-propttt}(e) shows  it is possible to compute $\jstar$ as the unique solution of the problem of maximizing $\sum_{i=1}^n \b_iJ(i)$ over all $J=\big(J(1),\ldots,J(n)\big)$ such that $J\le TJ$, where $\b_1,\ldots,\b_n$ are any positive scalars. This problem can be written as
\begin{align} \label{eq-lpproblem}
\hbox{\rm maximize}\quad &\sum_{i=1}^n \b_i J(i)\cr
\hbox{\rm subject to\ \ }
&J(i)\le b(i,u)
+ \sum_{j=1}^n A_{ij}(u)J(j),\notag\\
&\quad \quad \quad i =1,\ldots,n, \quad u\in U(i),
\end{align}
and it  is a linear program if each $U(i)$ is a finite set.

\section{Case of Finite Cost noncontractive Policies}

\pn We will now  eliminate Assumption \ref{assumptionthreeone}, thus allowing noncontractive policies with real-valued cost functions. We will prove results that are weaker yet useful and substantial.
An important notion in this regard is the optimal cost that can be achieved with contractive policies only, i.e., the vector  $\hat J$ with components given by
\begin{equation} \label{eq-optcontractive}
\hat J(i)=\inf_{\m:\, \hbox{contractive}}J_\m(i),\qquad i=1,\ldots,n.
\end{equation}
We will  show that $\hat J$ is a solution of Bellman's equation, while $\jstar$ need not be. 
To this end, we give an important property of noncontractive policies in the following proposition.

\begin{proposition}\label{prp-propnoncontractive}
If $\m$ is a noncontractive policy and all the components of $b_\m$ are strictly positive, then there exists at least one state $i$ such that the corresponding component of the vector
$\sum_{k=0}^{\infty} A_\m^k b_\m$
is $\infty$.
\end{proposition}

\proof According to the Perron-Frobenius Theorem, the nonnegative matrix $A_\m$ has a real eigenvalue $\l$, which is equal to its spectral radius, and an associated nonnegative eigenvector $\xi\ne0$ (see e.g., [22], Chapter 2, Prop.\ 6.6). Choose $\g>0$ to be such that $b_\m\ge \g\xi$, so that 
$$\sum_{k=0}^{\infty} A_\m^k b_\m\ge \g\sum_{k=0}^{\infty} A_\m^k \xi=\g\lf(\sum_{k=0}^{\infty} \l^k\ri) \xi.$$
Since some component of $\xi$ is positive while $\l\ge1$ (since $\m$ is noncontractive),  the corresponding component of the infinite sum on the right is infinite, and the same is true for the corresponding component of the vector
$\sum_{k=0}^{\infty} A_\m^k b_\m$ on the left. \qed

\subsection{The $\d$-Perturbed Problem}

We now introduce a perturbation line of analysis, also used in [1] and [26], whereby we add a constant $\d>0$ to all components of $b_\m$, thus obtaining what we call the {\it $\d$-perturbed affine monotonic model\/}. We denote by $J_{\m,\d}$ and $\jstar_\d$ the cost function of $\m$ and the optimal cost function of the $\d$-perturbed model, respectively. We have the following proposition.

\begin{proposition}\label{prp-propdeltaoptam}
Let Assumptions \ref{assumptiontwoone} and \ref{assumptiontwotwo} hold. Then for each $\d>0$:
\begin{itemize}
\item [(a)] $\jstar_\d$ is the unique solution within $\rn_+$ of the equation 
$$J(i)=(TJ)(i)+\d,\qquad i=1,\ldots,n.$$ 
\item [(b)] A policy $\m$ is optimal for the $\d$-perturbed problem (i.e., $J_{\m,\d}=\jstar_\d$) if and only if  $T_{\m}\jstar_\d =T\jstar_\d$. Moreover, for the $\d$-perturbed problem, all optimal policies are contractive and there exists at least one contractive policy that is optimal.
\item [(c)] The optimal cost function over contractive policies $\hat J$ [cf.\ Eq.\ \eqref{eq-optcontractive}] satisfies 
$$\hat J(i)=\lim_{\d\downarrow0}\jstar_{\d}(i),\qquad i=1,\ldots,n.$$ 
\item [(d)] If the control constraint set $U(i)$ is finite for all states $i=1,\ldots,n$, there exists a contractive policy $\hat \m$ that attains the minimum over all contractive policies, i.e., $J_{\hat \m}=\hat J$.
\end{itemize}
\end{proposition}

\proof (a), (b) By Prop.\ \ref{prp-propnoncontractive}, we have that Assumption \ref{assumptionthreeone} holds for the $\d$-perturbed problem. The results  follow by applying Prop.\ \ref{prp-propttt} [the equation of part (a) is Bellman's equation for the $\d$-perturbed problem]. 

\smskip
\pn (c) For an optimal contractive policy $\m^*_\d$ of the $\d$-perturbed problem [cf.\ part (b)], we have for all $\m'$ that are contractive
$$\hat J=\inf_{\m:\, \hbox{contractive}}J_\m\le J_{\m^*_\d}\le J_{\m^*_\d,\d}= \jstar_\d\le J_{\m',\d}.
$$
Since for every contractive policy $\m'$, we have
$\lim_{\d\downarrow0}J_{\m',\d}=J_{\m'},$
it follows that
$$\hat J\le \lim_{\d\downarrow0}J_{\m^*_\d}\le  J_{\m'},\qquad \forall\ \m': \hbox{contractive}.
$$
By taking the infimum over all $\m'$ that are contractive, the result follows. 

\smskip
\pn (d) Let $\{\d_k\}$ be a positive sequence with $\d_k\downarrow0$, and consider a corresponding sequence $\{\m_{k}\}$ of optimal contractive policies for the $\d_k$-perturbed problems. Since the set of contractive policies is finite, some policy $\hat \m$ will be repeated infinitely often within the sequence $\{\m_{k}\}$, and since $\{\jstar_{\d_k}\}$ is monotonically nonincreasing, we will have
$$\hat J\le J_{\hat \m}\le \jstar_{\d_k},$$
for all $k$ sufficiently large. Since by part (c), $\jstar_{\d_k}\downarrow\hat J$, it follows that $J_{\hat \m}=\hat J$.
 \qed

\subsection{Main Results}

\begin{figure}\vspace{-0pt}\hspace{1pc}
\scalebox{0.7}{\includegraphics[width=12cm]{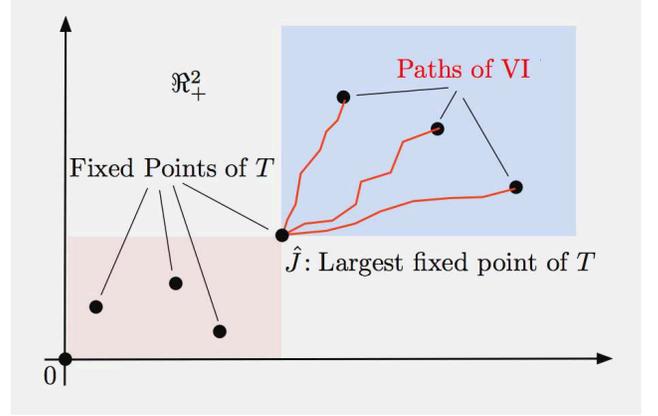}}\vspace{10pt}
\caption{Schematic illustration of  Prop.\ \ref{prp-propfixedviaff} for a problem with two states. The optimal cost function over contractive policies, $\skew6\hat J$, is the largest solution of Bellman's equation, while VI converges to $\skew6\hat J$ starting from $J\ge \skew6\hat J$.
}
\label{fig-Beleq}
\end{figure}

We now show that $\hat J$ is the largest fixed point of $T$ within $\rn_+$. This is the subject of the next proposition, which also provides a convergence result for VI as well as an optimality condition; see Fig.\ \ref{fig-Beleq}.

\begin{proposition}[Bellman's Equation,  Value Iteration, and Optimality Conditions]\label{prp-propfixedviaff}
Let Assumptions \ref{assumptiontwoone} and \ref{assumptiontwotwo} hold. Then:
\begin{itemize}
\item [(a)] The optimal cost function over contractive policies, $\hat J$, is the largest fixed point of $T$ within $\rn_+$, i.e., $\hat J$ is a fixed point that belongs to $\rn_+$, and if $J'\in\rn_+$ is another fixed point, we have $J'\le \hat J$. 
\item [(b)] We have $T^kJ\to \hat J$ for every
 $J\in \rn_+$ with $J\ge \hat J$.
\item [(c)] Let $\m$ be a contractive policy. Then  $\m$ is optimal within the class of contractive policies (i.e., $J_\m=\hat J$) if and only if $T_{\m} \hat J =T\hat J$. 
\end{itemize}
\end{proposition}

\proof (a), (b) For all contractive $\m$, we have
$J_\m=T_\m J_\m\ge T_\m\hat J\ge T\hat J.$
Taking the infimum over contractive $\m$, we obtain $\hat J\ge T\hat J$.
Conversely, for all $\d>0$ and $\m\in{\cal M}$, we have
$$\jstar_\d=T\jstar_\d+\d e\le T_\m \jstar_\d+\d e.$$
Taking limit as $\d\downarrow0$, and using  Prop.\ \ref{prp-propdeltaoptam}(c), we obtain $\hat J\le T_\m\hat J$ for all $\m\in{\cal M}$. Taking infimum over $\m\in{\cal M}$, it follows that $\hat J\le T\hat J$.  Thus $\hat J$ is a fixed point of $T$. 

For all $J\in\rn$ with $J\ge \hat J$ and contractive $\m$, we have by using the relation $\hat J=T\hat J$ just shown,
$$\hat J=\lim_{k\to\infty}T^k\hat J\le \lim_{k\to\infty}T^k J\le \lim_{k\to\infty}T_\m^k J=J_\m.$$
Taking the infimum over all contractive $\m$, we obtain
$$\hat J\le \lim_{k\to\infty}T^k J\le \hat J,\qquad \forall\ J\ge \hat J.$$
This proves that $T^kJ\to\hat J$ starting from
any $J\in \rn_+$ with $J\ge \hat J$. Finally, let $J'\in\rn_+$ be another fixed point of $T$, and let $J\in \rn_+$ be such that $J\ge \hat J$ and  $J\ge J'$. Then $T^kJ\to\hat J$, while $T^kJ\ge T^kJ'=J'$. It follows that $\hat J\ge J'$.

\smskip
\pn (c) If $\m$ is a contractive policy with $J_\m=\hat J$, we have
$\hat J=J_{\m}=T_{\m}J_{\m}=T_{\m}\hat J,$
so, using  also the relation $\hat J=T\hat J$ [cf.\ part (a)], we obtain $T_{\m}\hat J=T\hat J$.
 Conversely, if  $\m$  satisfies $T_{\m}\hat J=T\hat J$, then from part (a), we have $T_\m \hat J=\hat J$ and hence $\lim_{k\to\infty}T_\m^k\hat J = \hat J$. Since $\m$ is contractive, we obtain $J_\m=\lim_{k\to\infty}T_\m^k\hat J$, so  $J_\m=\hat J$.
\qed
\smskip

Note that it is possible that there exists a noncontractive policy $\m$ that is strictly suboptimal and yet satisfies the optimality condition $T_\m \jstar =T\jstar$ (there are simple deterministic shortest path examples with a zero length cycle that can be used to show this; see [1], Section 3.1.2). Thus contractiveness of $\m$ is an essential assumption in Prop.\ \ref{prp-propfixedviaff}(c). 

The following proposition shows that  starting from any $J\ge \hat J$, the convergence rate of VI to $\hat J$ is linear. The proposition also provides a corresponding error bound. The proof is very similar to a corresponding result of [26] and will not be given.

\begin{proposition}[Convergence Rate of VI]\label{prp-propvirateconvam}
Let Assumptions \ref{assumptiontwoone} and \ref{assumptiontwotwo} hold, and assume that there exists a contractive policy $\hat \m$ that is optimal within the class of contractive policies, i.e., $J_{\hat \m}=\hat J$. Then
$$\big\|TJ-\hat J\|_v\le \b\|J-\hat J\|_v,\qquad \forall\ J\ge \hat J,$$
where $\|\cdot\|_v$ is a weighted sup-norm for which $T_{\m^*}$ is a contraction and $\b$ is the corresponding modulus of contraction. Moreover, we have
$$\|J-\hat J\|_v\le {1\over 1-\b}\max_{i=1,\ldots,n}{J(i)-(TJ)(i)\over v(i)},\qquad \forall\ J\ge \hat J.$$\end{proposition}

We note that if $U(i)$ is infinite it is possible that $\hat J=\jstar$, but the only optimal policy is noncontractive, even if the compactness Assumption \ref{assumptiontwotwo} holds. This is shown in the following example, which is adapted from the paper [26] (Example 2.1).

\begin{figure}\vspace{-0pt}\hspace{2,5pc}
\scalebox{0.7}{\includegraphics[width=9cm]{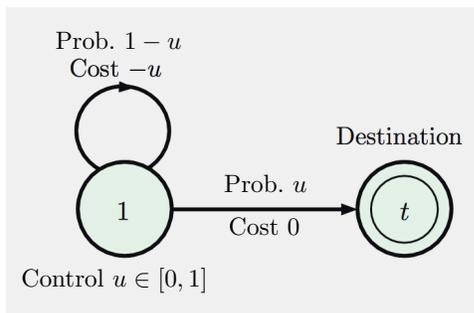}}\vspace{10pt}
\caption{The exponential cost SSP problem with a single state of Example \ref{examplecounterexp}.}
\label{fig-exp-cost}
\end{figure}

\begin{example}[A Counterexample on the Existence of an Optimal Contractive Policy] \label{examplecounterexp}
Consider an exponential cost SSP problem with a single state 1 in addition to the termination state $t$; cf.\ Fig.\ \ref{fig-exp-cost}. At state 1 we must choose $u\in[0,1]$. Then, we terminate at no cost [$g(1,u,t)=0$ in Eq.\ \eqref{eq-hexponspec}] with probability $u$, and we stay at state 1 at cost $-u$ [i.e., $g(1,u,1)=-u$ in Eq.\ \eqref{eq-hexponspec}] with probability $1-u$. We have $b(i,u)=u\exp{(0)}$ and $A_{11}(u)=(1-u)\exp{(-u)}$, so that 
$$H(1,u,J)=u+(1-u)\exp{(-u)}J.$$
Here there is a unique noncontractive policy $\m'$: it chooses $u=0$ at state 1, and has cost $J_{\m'}(1)=1$. Every policy $\m$ with $\m(1)\in (0,1]$ is  contractive, and $J_\m$ can be obtained by solving the equation $J_\m=T_\m J_\m$, i.e.,
$$J_\m(1)=\m(1)+\big(1-\m(1)\big)\exp{\big(-\m(1)\big)}J_\m(1).$$
We thus obtain 
$$J_\m(1)={\m(1)\over \big(1-\m(1)\big)\exp{\big(-\m(1)\big)}}.$$
It can be seen that $\skew5\hat J(1)=J^*(1)=0$, but there exists no optimal policy, and no optimal policy within the class of contractive policies. 
\end{example}

Let us also show that generally, under Assumptions \ref{assumptiontwoone} and \ref{assumptiontwotwo}, $\jstar$ need not be a fixed point of $T$. The following is a straightforward adaptation of Example 2.2 of [26].

\begin{figure}\vspace{-0pt}\hspace{0.3pc}
\scalebox{0.7}{\includegraphics[width=12cm]{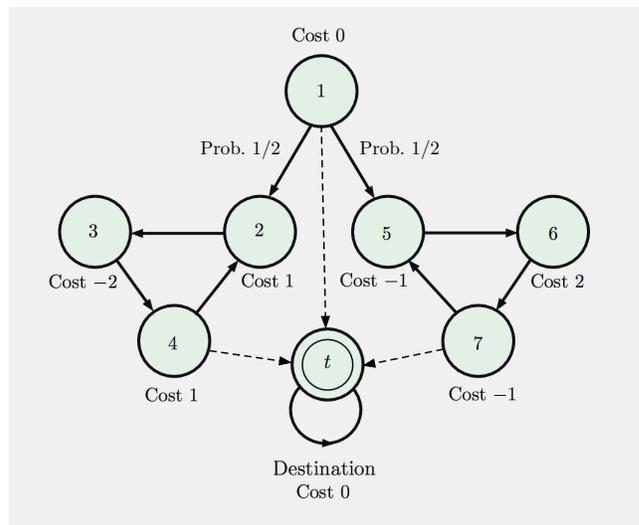}}\vspace{10pt}
\caption{An example of a noncontractive policy $\m$, where  $J_\m$ is not a fixed point of $T_\m$. All transitions under $\m$ are shown with solid lines. These transitions are deterministic, except at state 1 where the next state is 2 or 5 with equal probability $1/2$. There are additional transitions from nodes 1, 4, and 7 to the destination (shown with broken lines) with cost $c>2$, which create a suboptimal contractive policy. We have $J^*=J_\m$ and $J^*$ is not a fixed point of $T$.}
\label{fig-belcounterex1s}
\end{figure}

\begin{example}[An Exponential Cost SSP Problem Where  $\jstar$ is not a Fixed Point of $T$] \label{examplebelcounterexp}
Consider the exponential cost SSP problem of Fig.\ \ref{fig-belcounterex1s}, involving a noncontractive policy $\m$ whose transitions are marked by solid lines in the figure and form the two zero length cycles shown. All the transitions under $\m$ are deterministic, except at state 1 where the successor state is 2 or 5 with equal probability $1/2$. We assume that the cost of the policy for a given state is the expected value of the exponential of the finite horizon path length. We first calculate $J_\m(1)$. Let $g_k$ denote the cost incurred at time $k$, starting at state 1, and let $s_N(1)=\sum_{k=0}^{N-1} g_k$ denote the $N$-step accumulation of $g_k$ starting from state 1. We have
$$s_N(1)=0\qquad \hbox{if $N=1$ or $N=4+3t$, $t=0,1,\ldots$},$$
and 
\begin{align*}
s_N(1)=1\ \hbox{or $s_N(1)=-1$ with probability 1/2 each}\\
\qquad \qquad \hbox{if $N=2+3t$ or $N=3+3t$, $t=0,1,\ldots$.}
\end{align*}
Thus
$$J_\m(1)=\limsup_{N\to\infty}E\lf\{e^{s_N(1)}\ri\}={1\over 2} (e^1+e^{-1}).$$
On the other hand, a similar (but simpler) calculation shows that
$$J_\m(2)=J_\m(5)=e^1,$$
(the $N$-step accumulation of $g_k$ undergoes a cycle $\{1,-1,0,1,-1,0,\ldots\}$ as $N$ increases starting from state 2, and  undergoes a cycle $\{-1,1,0,-1,1,0,\ldots\}$ as $N$ increases starting from state 5).
Thus the Bellman equation at state 1, 
$$J_\m(1)={1\over 2} \big(J_\m(2)+J_\m(5)\big),$$
is not satisfied, and $J_\m$ is not a fixed point of $T_\m$. If for  $i=1,4,7,$ we have transitions (shown with broken lines) that lead from $i$ to $t$ with a cost $c>2$, the corresponding contractive policy is strictly suboptimal, so that $\m$ is optimal, but $J_\m=J^*$ is not a fixed point of $T$.
\end{example}

\subsection{Computational Methods}

\pn Regarding computational methods, Prop.\ \ref{prp-propfixedviaff}(b) establishes the validity of the VI algorithm that generates the sequence $\{T^kJ\}$ and converges to $\hat J$, starting from any initial $J\in \rn_+$ with $J\ge \hat J$. Moreover, Prop.\ \ref{prp-propvirateconvam} yields a linear rate of convergence result for this VI algorithm, assuming that there exists a contractive policy $\hat \m$ that is optimal within the class of contractive policies. Convergence to $\hat J$ starting from within the region $\{J\mid 0\le J\le \hat J\}$ cannot be guaranteed, since there may be fixed points other than $\hat J$ within that region. There are also PI algorithms that converge to $\hat J$. As an example, we note a PI algorithm with perturbations for abstract DP problems developed in Section 3.3.3 of [1], which can be readily adapted to affine monotonic problems. Finally, it is possible to compute $\hat J$ by solving a linear programming problem, in the case where the control space $U$ is finite, by using the following proposition.

\begin{proposition}\label{prp-proppertlp}
Let Assumptions  \ref{assumptiontwoone} and \ref{assumptiontwotwo} hold. Then if a vector $J\in \rn$  satisfies $J\le TJ$, it also satisfies $J\le \hat J $.
\end{proposition}

\proof Let $J\le TJ$ and $\d>0$. We have $J\le TJ+\d e=T_\d J$, and hence $J\le T^k_\d J$ for all $k$. Since the infinite cost conditions hold for the $\d$-perturbed problem, it follows that $T^k_\d J\to \jstar_\d$, so $J\le \jstar_\d$. By taking $\d\downarrow 0$ and using  Prop.\ \ref{prp-propdeltaoptam}(c), it follows that $J\le \hat J$. \qed
\smskip

The preceding proposition shows that $\hat J$ is the unique solution of the problem of maximizing $\sum_{i=1}^n \b_iJ(i)$ over all $J=\big(J(1),\ldots,J(n)\big)$ such that $J\le TJ$, where $\b_1,\ldots,\b_n$ are any positive scalars, i.e., the problem of Eq.\ \eqref{eq-lpproblem}. This problem is a linear program if each $U(i)$ is a finite set.

\vskip-1.5pc

\section{Concluding Remarks}

In this paper we have expanded the SSP methodology to affine monotonic models that are characterized by an affine mapping from the set of nonnegative functions  to itself. These models include among others, multiplicative and risk-averse exponentiated cost models. We have used the conceptual framework of semicontractive DP, based on the notion of a contractive policy, which generalizes the notion of a proper policy in SSP. We have provided extensions of the basic analytical and algorithmic results of SSP problems, and we have illustrated their exceptional behavior within our broader context.

Another case of affine monotonic model that we have not considered, is the one obtained when $\bar J\le0$ and 
$$b(i,u)\le0,\qquad A_{ij}(u)\ge0,\qquad \forall\ i,j=1,\ldots,n,\ u\in U(i),$$
so that $T_\m$ maps the space of nonpositive functions into itself. This case has different character from the case $\bar J\ge0$ and $b(i,u)\ge0$ of this paper, in analogy with the well-known differences in structure between stochastic optimal control problems with nonpositive and nonnegative cost per stage.

\vskip-1.5pc

\section{References}

\def\ref{\vskip1pt\pn}

\def\refer{\ref}

\ref[1] Bertsekas, D.\ P., 2013.\ Abstract Dynamic Programming, Athena Scientific, Belmont, MA.

\ref [2] Denardo, E.\ V., 1967.\  ``Contraction Mappings in the Theory Underlying
Dynamic Programming," SIAM Review, Vol.\ 9, pp.\ 165-177.

\ref[3] Denardo, E.\ V., and Mitten, L.\ G., 1967.\ 
``Elements of Sequential Decision Processes," J.\ Indust.\ Engrg., Vol.\ 18, pp.\ 106-112.

\ref [4] Bertsekas, D.\ P., 1975.\  ``Monotone Mappings in
Dynamic Programming," 1975 IEEE Conference on Decision and Control, pp.\
20-25.

\ref [5] Bertsekas, D.\ P., 1977.\  ``Monotone Mappings with Application in
Dynamic Programming," SIAM J.\ on Control and Opt., Vol.\ 15, pp.\
438-464.

\ref [6]  Bertsekas, D.\ P., and Shreve, S.\ E., 1978.\  Stochastic Optimal
Control:  The Discrete Time Case, Academic Press, N.\ Y.; may be downloaded 
from http://web.mit.edu/dimitrib/www/home.html

\ref [7] Blackwell, D., 1965.\  ``Positive Dynamic Programming," Proc.\ Fifth Berkeley Symposium Math.\ Statistics and Probability, pp.\ 415-418.

\ref [8] Strauch, R., 1966.\  ``Negative Dynamic Programming," Ann.\ Math.\
Statist., Vol.\ 37, pp.\ 871-890.

\ref[9]  Verdu, S., and Poor, H.\ V., 1984.\  ``Backward, Forward, and
Back\-ward-Forward Dynamic Programming Models under Commutativity Conditions,"
Proc.\ 1984 IEEE Decision and Control Conference,
Las Vegas, NE,  pp.\ 1081-1086.

\ref[10] Szepesvari, C., 1998.\ Static and Dynamic Aspects of Optimal Sequential Decision Making, Ph.D.\ Thesis, Bolyai Institute of Mathematics, Hungary.

\ref[11] Szepesvari, C., 1998.\ ``Non-Markovian Policies in Sequential Decision Problems,"
Acta Cybernetica, Vol.\ 13, pp.\ 305-318.

\ref[12] Bertsekas, D.\ P., 2012.\ Dynamic Programming and Optimal Control, Vol.\ II, 4th Edition: Approximate Dynamic Programming, Athena Scientific, Belmont, MA.

\ref[13] Bertsekas, D.\ P., and Yu, H., 2010.\ ``Asynchronous Distributed Policy Iteration in Dynamic Programming,"  Proc.\ of Allerton Conf.\ on Com., Control and Comp.,  Allerton Park, Ill, pp.\ 1368-1374.

\refer[14] Bertsekas, D.\ P., and Yu, H., 2012.\ ``Q-Learning and Enhanced Policy Iteration in Discounted Dynamic Programming,"  Math.\ of OR, Vol.\ 37, pp.\ 66-94.

\ref[15] Yu, H., and Bertsekas, D.\ P., 2013.\ ``Q-Learning and Policy Iteration Algorithms for Stochastic Shortest Path Problems," Annals of Operations Research, Vol.\ 208, pp.\ 95-132.

\ref [16] Eaton, J.\ H., and Zadeh, L.\ A., 1962.\  ``Optimal Pursuit Strategies
in Discrete State Probabilistic Systems," Trans.\ ASME Ser.\ D.\ J.\ Basic Eng.,
Vol.\ 84, pp.23-29.

\ref[17] Pallu de la Barriere, R., 1967.\ Optimal Control Theory, Saunders, Phila; republished by Dover, N. Y., 1980.

\ref [18] Derman, C., 1970.\ Finite State Markovian Decision Processes,
Academic Press, N.\ Y.

\ref [19] Whittle, P., 1982.\  Optimization Over Time, Wiley, N.\ Y., Vol.\
1, 1982, Vol.\ 2, 1983.

\ref[20] Kallenberg, L.\ C.\ M.\ 1983.\ Linear Programming and Finite Markovian Control Problems, MC Tracts 148, Amterdam.

\ref[21] Bertsekas, D.\ P., 1987.\ Dynamic Programming: Deterministic and Stochastic Models, Prentice-Hall, Englewood Cliffs, N.\ J.

\refer[22] Bertsekas, D.\ P., and Tsitsiklis, J.\ N., 1989.\ Parallel and
Distributed Computation: Numerical Methods, Prentice-Hall, Englewood Cliffs,
N.\ J.

\ref[23] Altman, E., 1999.\ Constrained Markov Decision Processes, CRC Press, Boca Raton, FL.

\ref[24] Hernandez-Lerma, O., and Lasserre, J.\ B., 1999.\ Further Topics on Discrete-Time Markov Control Processes, Springer, N.\ Y.

\ref [25]  Bertsekas, D.\ P., and Tsitsiklis, J.\ N., 1991.\ ``An Analysis of
Stochastic Shortest Path Problems,"
Math.\ of OR, Vol.\ 16, pp.\ 580-595.

\ref[26] Bertsekas, D.\ P., and Yu, H., 2016.\ ``Stochastic Shortest Path Problems Under Weak Conditions,"  Lab.\ for Information and Decision Systems Report LIDS-2909, August 2013, revised March 2015 and January 2016. 

\ref[27] Denardo, E.\ V., and Rothblum, U.\ G., 1979.\ ``Optimal Stopping, Exponential Utility, and Linear Programming," Math.\ Programming, Vol.\ 16, pp.\ 228-244.

\ref[28] Patek, S.\ D., 2001.\ ``On Terminating Markov Decision Processes with a Risk Averse Objective Function," Automatica, Vol.\ 37, pp.\ 1379-1386.

\ref [29] Puterman, M.\ L., 1994.\  Markov Decision Processes: Discrete Stochastic Dynamic Programming, J.\ Wiley, N.\ Y.

\ref[30] Rothblum, U.\ G., 1984.\ ``Multiplicative Markov Decision Chains," Math.\ of OR, Vol.\ 9, pp.\ 6-24.

\ref [31] Veinott, A.\ F., Jr., 1969.\  ``Discrete Dynamic Programming with
Sensitive Discount Optimality Criteria," Ann.\ Math.\ Statist., Vol.\ 40, pp.\
1635-1660.

\ref [32]  Bertsekas, D.\ P., and Tsitsiklis, J.\ N., 1996.\ Neuro-Dynamic
Programming, Athena Scientific, Belmont, MA.

\ref[33] Tseng, P., 1990.\ ``Solving $H$-Horizon, Stationary Markov Decision
Problems in Time Proportional to $\log (H)$,'' Operations Research Letters, Vol.\ 9,
pp.\ 287-297.

\ref[34] Littman, M.\ L., 1996.\ Algorithms for Sequential Decision Making, Ph.D.\ thesis, Brown University, Providence, R.\ I.

\ref[35]  Yu, H.,  and Bertsekas, D.\ P., 2013.\ ``A Mixed Value and Policy Iteration Method for Stochastic Control with Universally Measurable Policies," 
arXiv preprint arXiv:1308.3814, to appear in Math.\ of OR.

\begin{IEEEbiography}{Dimitri P.\ Bertsekas} studied engineering at the National Technical University of Athens, Greece, obtained his MS in electrical engineering at the George Washington University, Wash. DC in 1969, and his Ph.D. in system science in 1971 at the Massachusetts Institute of Technology.

Dr. Bertsekas has held faculty positions with the Engineering-Economic Systems Dept., Stanford University (1971-1974) and the Electrical Engineering Dept. of the University of Illinois, Urbana (1974-1979). Since 1979 he has been teaching at the Electrical Engineering and Computer Science Department of the Massachusetts Institute of Technology (M.I.T.), where he is currently McAfee Professor of Engineering. He consults regularly with private industry and has held editorial positions in several journals. His research at M.I.T. spans several fields, including optimization, control, large-scale computation, and data communication networks, and is closely tied to his teaching and book authoring activities. He has written numerous research papers, and fourteen books, several of which are used as textbooks in MIT classes.

Professor Bertsekas was awarded the INFORMS 1997 Prize for Research Excellence in the Interface Between Operations Research and Computer Science for his book ``Neuro-Dynamic Programming" (co-authored with John Tsitsiklis), the 2001 ACC John R. Ragazzini Education Award, the 2009 INFORMS Expository Writing Award, the 2014 ACC Richard E. Bellman Control Heritage Award, the 2014 Khachiyan Prize, and the SIAM/MOS 2015 George B. Dantzig Prize. In 2001, he was elected to the United States National Academy of Engineering.

Dr. Bertsekas' recent books are ``Dynamic Programming and Optimal Control: 4th Edition" (2017), ``Noninear Programming: 3rd Edition" (2016),  ``Convex Optimization Algorithms" (2015), and ``Abstract Dynamic Programming" (2013) all published by Athena Scientific. 
\end{IEEEbiography}

%%%%%%% CURRENT END OF EDITED PAPER CONTENT %%%%%%%%%%%%%

%%%%%%% END OF PAPER CONTENT %%%%%%%%%%%%%

\end{document}